\documentclass[12pt,reqno]{amsart}
\usepackage{amsmath,amssymb}
 \makeatletter
 \evensidemargin  \oddsidemargin
 \makeatother

\begin{document}
\thispagestyle{empty}
 \setcounter{page}{1}
\begin{center}
{\large\bf On approximate ternary $m$-derivations and $\sigma$-homomorphisms}\vskip.20in

{\Small \bf $^1$A. G. Ghazanfari
and Z. Alizadeh$^2$} \vskip.20in \vskip 4mm

{\footnotesize $^{1}$Department of Mathematics, Lorestan University,P.O. Box 465, Khoramabad, Iran

E-mail: ghazanfari.a@lu.ac.ir; zahra.alizade43@gmail.com}\
\end{center}

\vskip 5mm
 \noindent{\footnotesize{\bf Abstract:} In this paper we introduce ternary modules over ternary algebras and
 using fixed point methods, we prove the stability and
super-stability of  ternary additive, quadratic, cubic and
quartic derivations and $\sigma$-homomorphisms in such structures for the functional
equation
\begin{equation*}
\begin{split}
&\quad
f(ax+y)+f(ax-y)= a^{m-2}[f(x+y)+f(x-y)]\\&+2(a^2-1)[a^{m-2}f(x)+\frac{(m-2)(1-(m-2)^2)}{6}f(y)].
\end{split}
\end{equation*}
for each $m=1,2,3,4$.
 \vskip.10in

\noindent{\bf Keywords:} Ternary algebras; stability; approximation; derivations; homomorphisms.

  \newtheorem{df}{Definition}[section]
  \newtheorem{rk}[df]{Remark}
   \newtheorem{lem}[df]{Lemma}
   \newtheorem{thm}[df]{Theorem}
   \newtheorem{pro}[df]{Proposition}
   \newtheorem{cor}[df]{Corollary}
   \newtheorem{ex}[df]{Example}
   \newtheorem{re}[df]{Remark}

 \setcounter{section}{0}
 \numberwithin{equation}{section}

\vskip 5mm

\section {\bf{Introduction}}\vskip 2mm

We recall that a nonempty set $G$ is said to be a {\it ternary groupoid} provided there exists on $G$ a
ternary operation $[\cdot,\cdot, \cdot] : G \times G\times
G\rightarrow G$, which is
denoted by $(G, [\cdot, \cdot, \cdot]).$ The ternary groupoid $(G,
[\cdot, \cdot, \cdot])$ is said to be {\it commutative} if
$[x_1,x_2,x_3]=[x_{\sigma(1)},x_{\sigma(2)},x_{\sigma(3)}]$ for
all $x_1, x_2, x_3 \in G$ and all permutations $\sigma$ of $ \{1,
2, 3\}.$ If a binary operation $\circ$ is defined on $G$ such that
$[x, y, z] = (x\circ y)\circ z$ for all $x, y, z \in G$, then we
say that $[\cdot, \cdot, \cdot]$ is derived from $\circ$. We say
that $(G, [\cdot, \cdot, \cdot])$ is a {\it ternary semigroup} if
the operation $[\cdot, \cdot, \cdot]$ is associative, i.e., if
$[[x, y, z], u, v] = [x, [y, z, u], v] = [x, y, [z, u, v]]$ holds
for all $x, y, z, u, v\in G$ (see \cite{B-B-K}).

  A {\it ternary Banach algebra} is a complex Banach space $A$, equipped with a
ternary product $(x, y,z) \rightarrow [x, y,z]$ of $A^3$ into $A$,
which is $\mathbb{C}$-linear in the outer variables, conjugate $\mathbb{C}$-linear
in the middle variable and associative in the sense that $[x,
y,[z,w,v]] = [x,[w,z, y],v] = [[x, y,z],w,v]$ and satisfies $ \|
[x, y,z] \| \leq  \|x\| \cdot \|y\| \cdot \|z\| $.

Consider the functional equation $\Im_1(f ) = \Im_2(f )~~(\Im)$ in
a certain general setting. A function $g$ is an approximate
solution of $(\Im)$ if $\Im_1(g)$ and $\Im_2(g)$ are close in some
sense. The Ulam stability problem asks whether or not there exists
a true solution of $(\Im)$ near $g$. A functional equation is said
to be {\it superstable} if every approximate solution of the
equation is an exact solution of the functional equation. The
problem of stability of functional equations originated from a
question of Ulam \cite{Ul} concerning the stability of group
homomorphisms:\vskip 2mm

Let $(G_1, \ast)$ be a group and $(G_2, \star ,d)$ be a metric
group with the metric $d(\cdot, \cdot)$. Given $ \epsilon> 0$,
does there exist a $\delta (\epsilon) > 0$ such that, if a mapping
$h:G_1\longrightarrow G_2$ satisfies the inequality
$$d(h(x \ast
y),h(x) \star h(y)) <\delta
$$ for all $x,y\in G_1$, then there
exists a homomorphism $H:G_1\rightarrow G_2$ with
$d(h(x),H(x))<\epsilon$ for all $x\in G_1$?\vskip 2mm

 If the
answer is affirmative, we say that the equation of homomorphism
$H(x \ast y) = H(x) \star H(y)$ is {\it stable}. The concept of
stability for a functional equation arises when we replace the
functional equation by an inequality which acts as a perturbation
of the equation. Thus the stability question of functional
equations is that how do the solutions of the inequality differ
from those of the given functional equation?

In 1941, Hyers \cite{Hy} gave a first affirmative answer to the
question of Ulam for Banach spaces.\vskip 2mm

 Let $X$ and $Y$ be
Banach spaces. Assume that $f:X\longrightarrow Y$ satisfies
$$\|f(x + y) - f(x) - f(y) \|\leq \epsilon
$$ for all $x , y\in X$ and
some $\epsilon > 0$. Then there exists a unique additive mapping
$T:X\longrightarrow Y$ such that $\| f(x) - T(x) \| \leq \epsilon$
for all $x \in X$.\vskip 2mm

 A generalized version of the theorem
of Hyers for approximately additive mappings was given by Aoki
\cite{Aok} in 1950 (see also \cite{Bou}). In 1978, a generalized
solution for approximately linear mappings was given by Th.M.
Rassias \cite{R1}. He considered a mapping $f:X\to Y$ satisfying
the condition
$$\|f(x+y)-f(x)-f(y)\|\leq\epsilon(\|x\|^p+\|y\|^p)$$
for all $x,y\in X$, where $\epsilon\geq0$ and $0\leq p<1$. This
result was later extended to all $p\neq1$ and generalized by Gajda
\cite{G}, Rassias and Semrl \cite{R-S}, Isac and Rassias
\cite{I-R}.

The problem when
$p=1$ is not true. Counterexamples for the corresponding assertion
in the case $p=1$ were constructed by Gadja \cite{G}, Rassias and
Semrl \cite{R-S}.

On the other hand, Rassias \cite{JR1,JR2} considered the
Cauchy difference controlled by a product of different powers of
norm. Furthermore, a generalization of Rassias theorems was
obtained by G$\check{a}$vruta \cite{Gv1}, who replaced
$$\epsilon(\parallel x\parallel^p+\parallel y\parallel^p)$$ and
$\epsilon\|x\|^p\|y\|^p$ by a general control function
$\varphi(x,y)$. In 1949 and 1951, Bourgin \cite{Bo,Bo2} is the
first mathematician dealing with stability of (ring) homomorphism
$f(xy) = f(x)f(y)$. The topic of approximation of functional
equations on Banach algebras was studied by a number of
mathematicians (see \cite{E1}--\cite{PR}).

The functional equation
\begin{equation}\label{e--2}
f(x+y)+f(x-y)=2f(x)+2f(y)
\end{equation}
is related to a symmetric bi--additive mapping \cite{Az}. It
is natural that this equation is called a {\it quadratic
functional equation}. For more details about various results
concerning such problems, the readers refer to
\cite{E}--\cite{K-R}.

In 2002, Jun and Kim \cite{J-K} introduced the following cubic
functional equation
\begin{equation}\label{e13}
f(2x+y)+f(2x-y)=2f(x+y)+2f(x-y)+12f(x)
\end{equation}
and they established the general solution and the generalized
Hyers--Ulam--Rassias stability for the  functional equation
(\ref{e13}). Obviously, the mapping $f(x)=cx^3$ satisfies the
functional equation $(\ref{e13})$, which is called the {\it cubic
functional equation}. In 2005, Lee et al. \cite{Lee} considered
the following functional equation
\begin{equation}\label{e14}
f(2x+y)+f(2x-y)=4f(x+y)+4f(x-y)+24f(x)-6f(y)
\end{equation}
It is easy to see that the mapping $f(x)=dx^4$ is a solution of
the functional equation $(\ref{e14})$, which is called the {\it
quartic functional equation}.\vskip 4mm

\section {\bf {Preliminaries}}\vskip 2mm

In 2007, Park \cite{P-c} investigated the generalized stability of a quadratic mapping $f:A\longrightarrow B$,
 which is called a $C^*$-{\it ternary quadratic mapping}
if $f$ is a quadratic mapping satisfies
\begin{equation}\label{xy}
f([x,y,z])=[f(x),f(y),f(z)]
\end{equation}
 for all $x,y,z\in A.$
Let $A$ be an algebra. An additive mapping $f:A \rightarrow A$ is
called a ring derivation if $f(xy) = xf(y) + f(x)y$ holds for all
$x, y \in A.$ If, in addition, $f(\lambda x) = \lambda f(x)$ for
all $x \in A$ and all $\lambda\in\mathbb{F}$, then $f$ is called a
linear derivation, where $\mathbb{F}$ denotes the scalar field of
$A$. The stability result concerning derivations between operator
algebras was first obtained by $\check{S}$emrl \cite{Sem}. In
\cite{Bad}, Badora proved the stability of functional equation
$f(xy) = xf(y) + f(x)y$, where $f$ is a mapping on normed algebra
$A$ with unit.

Recently, shagholi et al.\cite{Sha} proved the stability of ternary quadratic derivations on ternary Banach algebras. Also Moslehian had
investigated the stability and superstability of ternary derivations on $C^*$-ternary rings \cite{MO}.

the monomial $f(x)=ax^m ~(x \in \mathbb{R})$ is a
solution of the functional equation (\ref{e12}) for each
$m=1,2,3,4.$
\begin{equation}\label{e12}
\begin{split}
&\quad
f(ax+y)+f(ax-y)= a^{m-2}[f(x+y)+f(x-y)]\\&+2(a^2-1)[a^{m-2}f(x)+\frac{(m-2)(1-(m-2)^2)}{6}f(y)].
\end{split}
\end{equation}

For $m=1,2,3,4$, the functional
equation (\ref{e12}) is equivalent to the additive, quadratic, cubic and quartic
functional equation, respectively.

The general solution of the functional equation (\ref{e12}) for any fixed integers a with $a\neq 0,\pm1$, was obtained by Eshaghi Gordji et al.\cite{Es-A-Kh}. \vskip
2mm

In this paper, we study the
further generalized stability of ternary additive,
quadratic, cubic and quartic derivations and $\sigma$-homomorphisms over ternary Banach algebras
via fixed point method for the functional equation (\ref{e12}).
 Moreover, we establish the super-stability of this
functional equation by suitable control functions.\vskip 2mm

\begin{df}\label{df0} {\rm Let $A$ be a ternary algebra and $X$ be a vector space.\vskip 4mm
(i) $X$ is a left ternary $A$-module if mapping $ A\times A\times X \to X$ satisfies \vskip 2mm

(LTM 1) For each fixed $a\in A$, the mapping $x\to [a,b,x]$ is linear on $X$;

(LTM 2) For each fixed $x\in X$, the mapping $(a,b)\to [a,b,x]$ is bilinear on $A\times A$;

(LTM 3) For all $x\in X$ and $a,b,c,d\in A$, $[a, b, [c, d, x]]=[[a, b, c], d, x] = [a, [b, c, d], x] $.\vskip 4mm

 (ii) $X$ is a middle ternary $A$-module if mapping $ A\times X \times A \to X$ satisfies \vskip 2mm

 (MTM 1) For each fixed $a\in A$, the mapping $x\to [a,x,b]$ is linear on $X$;

 (MTM 2) For each fixed $x\in X$, the mapping $(a,b)\to [a,x,b]$ is bilinear on $A\times A$;

 (MTM 3) For all $x\in X$ and $a,b,c,d,e,f\in A$, $[a, [b, [c, x, d],e],f]=[[a, b, c], x,[d, e, f]] $.\vskip 4mm

(iii) $X$ is a right ternary $A$-module if mapping $ X\times  A\times A \to X$ satisfies \vskip 2mm

 (RTM 1) For each fixed $a\in A$, the mapping $x\to [x,a,b]$ is linear on $X$;

 (RTM 2) For each fixed $x\in X$, the mapping $(a,b)\to [x,a,b]$ is bilinear on $A\times A$;

 (RTM 3) For all $x\in X$ and $x,y,u,v\in A$,  $[[x,a,b],c,d]=[x, [a,b,c],d] = [x, a, [b,c,d]]$.\vskip 4mm

 (iv)  $X$ is called ternary $A$-module if $X$ is left ternary $A$-module, middle ternary $A$-module and right ternary $A$-module and if
       satisfies in the following condition:\vskip 2mm

  (TM) For any $a,b,c,d,e\in(A\cup X)$, which one of them is in $X$ and the rest are in $A$, $[[a, b, c], d, e] = [a,[b,c,d],e] = [a,b,[c,d,e]]$.
}
\end{df}\vskip 2mm

\begin{df}\label{df0} {\rm Let $A$ be a normed ternary algebra and $X$ be a vector space.\vskip 4mm
$X$ is said to be a normed left ternary $A$-module if $X$ is a left ternary $A$-module and also satisfies
in the following axiom: \vskip 2mm

(NLTM) $ \|[a,b,x] \| \leq  \|a\| \cdot \|b\| \cdot \|x\| $ for all $a,b\in A$ and $x\in X$.\vskip 4mm

Similarly normed middle ternary $A$-module, normed right ternary $A$-module and normed ternary $A$-module
are defined. A normed ternary $A$-module is called Banach ternary $A$-module if it is complete as a normed linear space.
}
\end{df}\vskip 2mm

\begin{df}\label{df1} {\rm Let $A,B$ are ternary algebra, $H:A\rightarrow B$ a function and $\sigma$ a permutation of $\{1,2,3\}$.
$H$ is said to be a ternary $\sigma$-hommomorphism if for all $a_1,a_2,a_3 \in A$
}
\begin{equation}\label{xy}
H([a_1,a_2,a_3])=[H(a_{\sigma(1)}),H(a_{\sigma(2)}),H(a_{\sigma(3)})]
\end{equation}
\end{df}
%%%%%%%%%%%%%%%%%%%%%%%%%%%%%%%%%%%%%%%%%%%%%%%%%%%%%%%%%%%%%%%%%%%%%%%%%%%%%%%%%%%%%%%%%%%%%%%%%%%%%%%%%%%%%%%%%%%%%
\begin{df}\label{df1} {\rm Let $A$ and $B$ be two ternary algebras.

(1) A mapping $f : A \to B$ is called a {\it ternary
additive $\sigma$-homomorphism} $($briefly, ternary
$1$-$\sigma$-homomorphism$)$ if $f$ is an additive mapping satisfying
$(\ref{xy})$  for all $a_1,a_2,a_3\in A$.\vskip 4mm

(2) A mapping $f : A \to B$ is called a {\it ternary
quadratic $\sigma$-homomorphism} $($briefly, ternary $2$-$\sigma$-homomorphism$)$
if $f$ is a quadratic mapping  satisfying $(\ref{xy})$  for all
$a_1,a_2,a_3\in A$.\vskip 4mm

(3) A mapping $f : A \to B$ is called a {\it ternary cubic
$\sigma$-homomorphism} $($briefly, ternary $3$-$\sigma$-homomorphism$)$ if $f$ is a
cubic mapping satisfying $(\ref{xy})$  for all $a_1,a_2,a_3\in A$.\vskip 4mm

(4) A mapping $f : A \to B$ is called a {\it ternary quartic
$\sigma$-homomorphism} $($briefly, ternary $4$-$\sigma$-homomorphism$)$ if $f$
is a quartic mapping satisfying $(\ref{xy})$  for all $a_1,a_2,a_3\in
A$.}
\end{df}\vskip 2mm

%%%%%%%%%%%%%%%%%%%%%%%%%%%%%%%%%%%%%%%%%%%%%%%%%%%%%%%%%%%%%%%%%%%%%%%%%%%%%%%%%%%%%%%%%%%%%%%%%%%%%%%%%%%%%%%%%%%%%%%%%%%%
\begin{df}\label{df1} {\rm Let $A$ be a ternary algebra and let $X$ be a Banach ternary $A$-module.\vskip 4mm

(1) A mapping $f : A \to X$ is called a {\it ternary
additive derivation} $($briefly, ternary
$1$-derivation$)$ if $f$ is an additive mapping that satisfies

$f([x,y,z])= [f(x),y,z]+[x,f(y),z]+[x,y,f(z)]$  for all $x,y,z\in A$.\vskip 4mm

(2) A mapping $f : A \to X$ is called a {\it ternary
quadratic derivation} $($briefly, ternary $2$-derivation$)$
if $f$ is a quadratic mapping that satisfies

$f([x,y,z])= [f(x),y^2,z^2]+[x^2,f(y),z^2]+[x^2,y^2,f(z)]$  for all
$x,y,z\in A$.\vskip 4mm

(3) A mapping $f : A \to X$ is called a {\it ternary cubic
derivation} $($briefly, ternary $3$-derivation$)$ if $f$ is a
cubic mapping that satisfies

$f([x,y,z])= [f(x),y^3,z^3]+[x^3,f(y),z^3]+[x^3,y^3,f(z)]$  for all $x,y,z\in A$.\vskip 4mm

(4) A mapping $f : A \to X$ is called a {\it ternary quartic
derivation} $($briefly, ternary $4$-derivation$)$ if $f$
is a quartic mapping that satisfies

$f([x,y,z])= [f(x),y^4,z^4]+[x^4,f(y),z^4]+[x^4,y^4,f(z)]$  for all $x,y,z\in
A$.}
\end{df}\vskip 2mm

 Now, we  state  the following notion of fixed point
theorem. For the proof, refer to Chapter 5 in
\cite{Rus} and \cite{Hy2,Ze}.  In 2003, C\v{a}dariu and Radu \cite{Cad1} proposed
a new method for obtaining the existence of exact solutions and
error estimations, based on the fixed point alternative (see also
\cite{C-R2004}--\cite{1}).\vskip 2mm

 Let $(X, d)$ be a generalized metric space.
We say that a mapping $T:X\rightarrow X$ satisfies a Lipschitz
condition if there exists a constant $L \geq 0$ such that $
d(Tx,Ty)\leq L d(x,y)$
 for all $x,y \in X,$ where the number $L$ is called the Lipschitz constant.
 If the
Lipschitz constant $L$ is less than $1$, then the mapping $T$ is
called a {\it strictly contractive mapping}. Note that the
distinction between the generalized metric and the usual metric is
that the range of the former is permitted to include the
infinity.\vskip 2mm

The following theorem is famous to fixed point theorm.\vskip 2mm

\begin{thm}\label{t2.1}
Suppose that $(\Omega,d)$ is a complete generalized metric space
 and $T:\Omega\rightarrow\Omega$ is a strictly contractive mapping
 with the Lipschitz constant $L$. Then, for
any $x\in\Omega$, either
$$
d(T^m x, T^{m+1} x)=\infty,\quad\forall m\geq0,
$$
 or there exists a natural number $m_{0}$ such that

$(1)$ $d(T^m x, T^{m+1} x)<\infty ~$ for all $m \geq m_{0}$;

$(2)$ the sequence $\{T^m x\}$ is convergent to a fixed point
$y^*$ of $~T$;

$(3)$ $y^*$ is the unique fixed point of $~T$ in
$\Lambda=\{y\in\Omega:d(T^{m_{0}} x, y)<\infty\};$

$(4)$ $d(y,y^*)\leq\frac{1}{1-L}d(y, Ty)$ for all $y\in\Lambda.$
\end{thm}
\vskip 4mm

\section {\bf{ Approximation of ternary $m$-derivation between ternary algebras
}}\vskip 2mm

In this section, we investigate the generalized  stability of
ternary $m$-derivation between ternary Banach algebras
for the functional equation $(\ref{e12})$.

Throughout this section, we suppose that $A$ is
ternary Banach algebra and $X$ is Banach ternary $A$-module. For convenience, we use the following
abbreviation: for any function $f:A\rightarrow X$,
$$
\aligned
\Delta_m{f}(x,y)&=f(ax+y)+f(ax-y)-a^{m-2}[f(x+y)+f(x-y)]\\
&\quad-2(a^2-1)[a^{m-2}f(x)+\frac{(m-2)(1-(m-2)^2)}{6}f(y)] \endaligned
$$ for all $x,y \in
X$ and any fixed integers $a$ with
$a\neq 0,\pm1$.\vskip 2mm

From now on, let $m$ be a positive integer less than $5$. \vskip
2mm

\begin{thm}\label{t3.1} Let $f:A\rightarrow X$ be a mapping for which there exist functions $\varphi_m:A\times A \rightarrow  [0,\infty)$
and $\psi_m:A\times A\times A \rightarrow  [0,\infty)$ such that
\begin{equation}\label{e31}
\|\Delta_m{f(x,y)}\|\leq \varphi_m(x,y),
\end{equation}
\begin{equation}\label{e32}
\|f([x,y,z])-[f(x),y^m,z^m]-[x^m,f(y),z^m]-[x^m,y^m,f(z)]\|\leq\psi_m(x,y,z)
\end{equation}
for all $x,y,z \in A.$ If there exists a constant $0<L<1$ such that
\begin{equation}\label{e33}
\varphi_m\Big(\frac{x}{a},\frac{y}{a}\Big)\leq
\frac{L}{|a|^m}\varphi_m(x,y),\end{equation}
\begin{equation}\label{e34}
\psi_m\Big(\frac{x}{a},\frac{y}{a},\frac{z}{a}\Big)\leq
\frac{L}{|a|
^{3m}}\psi_m(x,y,z)
\end{equation}
for all $x,y,z \in A,$ then there exists a unique ternary
$m$-derivation $\mathfrak{F}:A\rightarrow X$ such that
\begin{equation}\label{e35}
\|f(x)-\mathfrak{F}(x)\| \leq \frac{L}{2|a|^m(1-L)}\varphi_m(x,0)
\end{equation}
for all $x \in A.$
\end{thm}\vskip 2mm

\begin{proof}
It follows from $(\ref{e33})$ and $(\ref{e34})$ that
\begin{equation}\label{e36}
\lim_{n\rightarrow\infty}|a|^{mn}\varphi_m\Big(\frac{x}{a^n},\frac{y}{a^n}\Big)=0,
\end{equation}
\begin{equation}\label{e37}
\lim_{n\rightarrow\infty}|a|^{3mn}\psi_m\Big(\frac{x}{a^n},\frac{y}{a^n},\frac{z}{z^n}\Big)=0
\end{equation}for
all $x,y,z \in X$. By $(\ref{e36})$,
$\lim_{n\rightarrow\infty}|a|^{mn}\varphi_m(0,0)=0$ and so
$\varphi_m(0,0)=0$. Letting $x=y=0$ in $(\ref{e31})$, we get
$\|f(0)\| \leq \varphi_m(0,0)=0$ and so $f(0)=0$.

Let $\Omega=\{g:g:A \rightarrow X,~ g(0)=0\}$. We introduce a
generalized metric on $\Omega$ as follows:
$$d(g,h)=d_{\varphi_m}(g,h)=\inf\{K \in (0,\infty):~\|g(x)-h(x)\| \leq K \varphi_m(x,0),~\forall x \in A\}. \eqno \hspace{2cm}$$
It is easy to show that $(\Omega,d)$ is a generalized complete
metric space \cite{C-R2004}.

Now, we consider the mapping $T:\Omega\rightarrow\Omega~$ defined
by $~Tg(x)=a^m~~g(\frac{x}{a})$ for all $x \in A$ and $g \in
\Omega.$ Note that, for all $~g,h\in \Omega$ and $x\in A$,
\begin{align*}
d (g,h)< K &~\Rightarrow~~\|g(x)-h(x)\| \leq K\varphi_m(x,0)
\\&~\Rightarrow~~\|a^{m}g(\frac{x}{a})-a^{m}h(\frac{x}{a})\|\leq
|a|^{m}~K~\varphi_m(\frac{x}{a},0)
\\&~\Rightarrow~~\|a^{m}g(\frac{x}{a})-a^{m}h(\frac{x}{a})\|\leq
L~K~\varphi_m(x,0)
\\&~\Rightarrow~~d(Tg,Th)\leq L~K.
 \end{align*}
Hence we see that
$$d(Tg,Th)\leq L~d(g,h)$$
for all $g, h\in \Omega,$ that is, $T$ is a strictly self-mapping
of $\Omega$ with the Lipschitz constant $L$. Putting $y=0$ in
$(\ref{e31})$, we have
\begin{equation}\label{e38}
\|2f(ax)-2a^mf(x)\| \leq \varphi_m(x,0)
\end{equation}
for all $x \in A$ and so
\begin{equation*}
 \Big\|f(x)-a^mf\Big(\frac{x}{a}\Big)\Big\|\leq\frac{1}{2}\varphi_m\Big(\frac{x}{a},0\Big)
\leq\frac{L}{2|a|^m}~\varphi_m(x,0)
\end{equation*}
for all $x\in A,$ that is, $d(f,Tf)\leq\frac{L}{2|a|^m}<\infty$.

Now, from Theorem 2.5, it follows that there exists a fixed point
$\mathfrak{F}$ of $T$ in $\Omega$ such that
\begin{equation}\label{e39}
\mathfrak{F}(x)=\lim_{n\rightarrow\infty}a^{mn}f\Big(\frac{x}{a^{n}}\Big)
\end{equation}
 for all $x\in A$ since $\lim_{n\rightarrow\infty}d(T^nf,\mathfrak{F})=0.$

On the other hand, it follows from $(\ref{e31})$, $(\ref{e36})$
and $(\ref{e39})$ that
\begin{equation*}
\|\Delta_m \mathfrak{F}(x,y)\|=
\lim_{n\rightarrow\infty}|a|^{mn}\Big\|\Delta_m
f\Big(\frac{x}{a^{n}},\frac{y}{a^{n}}\Big)\Big\|\leq
\lim_{n\rightarrow\infty}|a|^{mn}\varphi_m\Big(\frac{x}{a^{n}},\frac{y}{a^{n}}\Big)=0
\end{equation*}
for all $x,y\in A$ and so $\Delta_m \mathfrak{F}(x,y)=0$.
 By the result in \cite{Es-A-Kh}, $\mathfrak{F}$ is $m-$mapping and so it follows from the
definition of $\mathfrak{F}$, $(\ref{e32})$ and $(\ref{e37})$ that
\begin{align*}
&\quad\|\mathfrak{F}([x,y,z])-[\mathfrak{F}(x),y^m,z^m]-[x^m,\mathfrak{F}(y),z^m]-[x^m,y^m,\mathfrak{F}(z))]\|\\
&=\lim_{n\rightarrow\infty}|a|^{3mn}\Big\|f\Big(\frac{[x,y,z]}{a^{3n}}\Big)-\Big[f\Big(\frac{x}{a^n}\Big),
\Big(\frac{y}{a^n}\Big)^m,\Big(\frac{z}{a^n}\Big)^m\Big]-\Big[\Big(\frac{x}{a^n}\Big)^m,f\Big(\frac{y}{a^n}\Big)
,\Big(\frac{z}{a^n}\Big)^m\Big]\\&~-\Big[\Big(\frac{x}{a^n}\Big)^m,\Big(\frac{y}{a^n}\Big)^m,f\Big(\frac{z}{a^n}\Big)\Big]\Big\|
\leq\lim_{n\rightarrow\infty}|a|^{3mn}\psi_m\Big(\frac{x}{a^n},\frac{y}{a^n},\frac{z}{a^n}\Big)=0
\end{align*}
 for all $x,y,z\in A$ and so

  $\mathfrak{F}([x,y,z])=[\mathfrak{F}(x),y^m,z^m]+[x^m,\mathfrak{F}(y),z^m]+[x^m,y^m,\mathfrak{F}(z))]$.

According to Theorem 2.2, since $\mathfrak{F}$ is the unique fixed
point of $T$ in the set $\Lambda=\{g\in\Omega:d(f,g)<\infty\},$
$\mathfrak{F}$ is the unique mapping such that
\begin{equation*}
 \|f(x)-\mathfrak{F}(x)\|\leq K~\varphi_m(x,0)
\end{equation*}
for all $x\in A$ and $K>0$. Again, using Theorem 2.2, we have
\begin{equation*}
d(f,\mathfrak{F}) \leq \frac{1}{1-L}d(f,Tf)\leq
\frac{L}{2|a|^m(1-L)}
\end{equation*}
and so
\begin{equation*}
\|f(x)-\mathfrak{F}(x)\| \leq\frac{L}{2|a|^m(1-L)}~\varphi_m(x,0)
\end{equation*}
for all $x \in A.$  This completes the proof.
\end{proof}
\vskip 2mm
%%%%%%%%%%%%%%%%%%%%%%%%%%%%%%%%%%%%%%%%%%%%%%%%%%%%%%%%%%%%%%%%%%%%%%%%%%%%%%%%%%%%%%%%%%%%%%%%%%%%%%%%
\begin{cor}\label{c3.2} Let $\theta,r,p$ be non--negative
real numbers with $r,p>m$ and $\frac{3p-r}{2}\geq m$. Suppose that
$f:A\rightarrow X$ is a mapping such that
\begin{equation}\label{b39}
\|\Delta_m f(x,y)\|\leq \theta(\|x\|^r+\|y\|^r),
\end{equation}
\begin{equation}\label{b39+}
\|f([x,y,z])-[f(x),y^m,z^m]-[x^m,f(y),z^m]-[x^m,y^m,f(z)]\|\leq
\theta(\|x\|^p.\|y\|^p.\|z\|^p)
\end{equation}
for all $x,y,z\in A.$ Then there exists a unique ternary
$m$-derivation $\mathfrak{F}:A\rightarrow X$ satisfying
\begin{equation*} \|f(x)-\mathfrak{F}(x)\|
\leq\frac{\theta}{2(|a|^r-|a|^m)}\|x\|^r
\end{equation*}
for all $x \in A.$
\end{cor}\vskip 2mm

\begin{proof}
The proof follows from Theorem \ref{t3.1} by taking
\begin{equation*}
\varphi_m(x,y):=\theta(\|x\|^r+\|y\|^r),\quad
\psi_m(x,y,z):=\theta(\|x\|^p\cdot\|y\|^p\cdot\|z\|^p)
\end{equation*}
for all $x,y,z \in A$. Then we can choose $L=|a|^{m-r}$ and so the
desired conclusion follows.
\end{proof}\vskip 2mm
%%%%%%%%%%%%%%%%%%%%%%%%%%%%%%%%%%%%%%%%%%%%%%%%%%%%%%%%%%%%%%%%%%%%%%%%%%%%%%%%%%%%%%%%%%%%%%%%%%%%%%%%%%%%%%%%%%%%%%%

\begin{rk}\label{rk}
{\rm Let  $f : A \to X $ be a mapping with $f(0)=0$ such that there
exist functions $\varphi_m:A\times A\rightarrow  [0,\infty)$ and
$\psi_m:A\times A\times A \rightarrow  [0,\infty)$ satisfying
(\ref{e31}) and $(\ref{e32})$. Let $0<L<1$ be a constant such that
$$\varphi_m(ax,ay)\leq |a|^mL\varphi_m(x,y),\quad
\psi_m(ax,ay,az)\leq |a|^{3m}L\psi_m(x,y,z)
$$
for all $x,y,z \in A.$ By the similar method as in the proof of
Theorem \ref{t3.1}, one can show that there exists a unique
ternary $m$-derivation $\mathfrak{F}:A\rightarrow X$
satisfying
\[\|f(x)-\mathfrak{F}(x)\| \leq
\frac{1}{2|a|^m(1-L)}\varphi_m(x,0)\] for all $x\in A.$
 For the case
 $$\varphi_m(x,y):=\delta+\theta(\|x\|^r+\|y\|^r),\quad
 \psi_m(x,y,z):=\delta+\theta(\|x\|^p\cdot\|y\|^p\cdot\|z\|^p),
 $$
  where $\theta, \delta$ are non--negative real
 numbers
 and $0<r,p<m$ and $\frac{3p-r}{2}\leq m,$
there exists a unique ternary $m$-derivation
$\mathfrak{F}:A\rightarrow X$ satisfying
\begin{equation*}
\|f(x)-\mathfrak{F}(x)\|
\leq\frac{\delta}{2(|a|^m-|a|^r)}+\frac{\theta}{2(|a|^m-|a|^r)}\|x\|^r
\end{equation*}
 for all $x \in A.$}
\end{rk}\vskip 2mm
%%%%%%%%%%%%%%%%%%%%%%%%%%%%%%%%%%%%%%%%%%%%%%%%%%%%%%%%%%%%%%%%%%%%%%%%%%%%%%%%%%%%%%%%%%%%%%%%%%%%%%%

In the following, we formulate and prove a theorem in
super-stability of ternary $m$-derivation  in ternary
algebras for the functional equation $(\ref{e12}).$\vskip 2mm

\begin{thm}\label{t3.5}
Suppose that there  exist functions $\varphi_m:A\times A
\rightarrow [0,\infty)$, $\psi_m:A\times A\times A \rightarrow
[0,\infty)$ and  a constant $0<L<1$ such that
\begin{equation}\label{Se3}
\varphi_m\Big(0,\frac{y}{a}\Big)\leq
\frac{L}{|a|^m}\varphi_m(0,y),\end{equation}
\begin{equation}\label{Se4}
\psi_m\Big(\frac{x}{a},\frac{y}{a},\frac{z}{a}\Big)\leq
\frac{L}{|a|^{3m}}\psi_m(x,y,z)
\end{equation}
for all $x,y,z \in A.$ Moreover, if $f : A \to X$ is a mapping
such that
\begin{equation}\label{Se5}
\|\Delta_m{f(x,y)}\|\leq \varphi_m(0,y),
\end{equation}
\begin{equation}\label{Se6}
\|f([x,y,z])-[f(x),y^m,z^m]-[x^m,f(y),z^m]-[x^m,y^m,f(z)]\|\leq\psi_m(x,y,z)
\end{equation}
for all $x,y,z \in A,$ then $f$ is a ternary
$m$-derivation.
\end{thm}\vskip 2mm

\begin{proof} It follows from $(\ref{Se3})$ and $(\ref{Se4})$ that
\begin{equation}\label{Se7}
\lim_{n\rightarrow\infty}|a|^{mn}\varphi_m\Big(0,\frac{y}{a^n}\Big)=0,
\end{equation}
\begin{equation}\label{Se8}
\lim_{n\rightarrow\infty}|a|^{3mn}\psi_m\Big(\frac{x}{a^n},\frac{y}{a^n},\frac{z}{a^n}\Big)=0
\end{equation}for
all $x,y,z \in A$. We have $f(0)=0$ since $\varphi_m(0,0)=0$.
Letting $y = 0$ in $(\ref{Se5})$, we get $f(ax) = a^{m}f(x)$ for
all $x \in A$. By using induction, we obtain $$f(a^nx) =
a^{mn}f(x)$$ for all $x \in A$ and $n \in \mathbb{N}$ and so
\begin{equation}\label{Se9}
f(x) = a^{mn}f\Big(\frac{x}{a^n}\Big)
\end{equation}
 for all $x \in A$ and $n \in
\mathbb{N}$. It follows from $(\ref{Se6})$ and $(\ref{Se9})$ that
\begin{equation}\label{Se10}
\begin{split}
&\quad\|f([x,y,z])-[f(x),y^m,z^m]-[x^m,f(y),z^m]-[x^m,y^m,f(z)]\|\\
&=
|a|^{3mn}\Big\|f\Big(\frac{[x,y,z]}{a^{3n}}\Big)-\Big[f\Big(\frac{x}{a^n}\Big),\Big(\frac{y}{a^n}\Big)^m,\Big(\frac{z}{a^n}\Big)^m\Big]-\Big[\Big(\frac{x}{a^n}\Big)^m,f\Big(\frac{y}{a^n}\Big)
,\Big(\frac{z}{a^n}\Big)^m\Big]\\&~-\Big[\Big(\frac{x}{a^n}\Big)^m,\Big(\frac{y}{a^n}\Big)^m,f\Big(\frac{z}{a^n}\Big)\Big]\Big\|
\leq |a|^{3mn}
\psi_m\Big(\frac{x}{a^n},\frac{y}{a^n},\frac{z}{a^n}\Big)
\end{split}
\end{equation}
 for all $x,y,z \in A$ and $n \in
\mathbb{N}$. Hence, letting $n\to\infty$ in $(\ref{Se10})$ and
using $(\ref{Se8})$, we have $f([x,y,z])=[f(x),y^m,z^m]+[x^m,f(y),z^m]+[x^m,y^m,f(z)]$
 for all $x,y,z \in A$.

On the other hand, we have
\begin{equation}\label{Se11}
\|\Delta_m f(x,y)\|=|a|^{mn} \Big\|\Delta_m
f\Big(\frac{x}{a^{n}},\frac{y}{a^{n}}\Big)\Big\|\leq
|a|^{mn}\varphi_m\Big(0,\frac{y}{a^{n}}\Big)
\end{equation}
 for all $x,y \in A$ and $n \in
\mathbb{N}$. Thus, letting $n\to\infty$ in $(\ref{Se11})$ and
using $(\ref{Se7})$, we have $\Delta_m f(x,y)=0$ for all $x,y \in
A$. Therefore, $f$ is a ternary $m$-derivation. This
completes the proof.
\end{proof}\vskip 2mm
%%%%%%%%%%%%%%%%%%%%%%%%%%%%%%%%%%%%%%%%%%%%%%%%%%%%%%%%%%%%%%%%%%%%%%%%%%%%%%%%%%%%%%%%%%%%%%%%%%%%%%%%%%%%%%%%%%%%%%%

\begin{cor}\label{c3.6} Let $\theta,r,s$ be non--negative
real numbers with $r>m$ and $s>3m$. If $f:A\rightarrow X$ is a
function such that
$$\|\Delta_m f(x,y)\|\leq \theta\|y\|^r,
$$
$$
\|f([x,y,z])-[f(x),y^m,z^m]-[x^m,f(y),z^m]-[x^m,y^m,f(z)]\|\leq
\theta(\|x\|^s+\|y\|^s+\|z\|^s)
$$ for all $x,y,z\in A,$ then $f$
is a ternary $m$-derivation.
\end{cor}\vskip 2mm
%%%%%%%%%%%%%%%%%%%%%%%%%%%%%%%%%%%%%%%%%%%%%%%%%%%%%%%%%%%%%%%%%%%%%%%%%%%%%%%%%%%%%%%%%%%%%%%%%%%%%%%%%%%%%%%%%%%%%%%

\begin{rk}\label{rk2} {\rm Let $\theta,r$ be non--negative real
numbers with $r<m$. Suppose that there exists a function
$\psi_m:A\times A\times A\rightarrow[0,\infty)$ and a constant
$0<L<1$  such that
$$\psi_m(ax,ay,az)\leq |a|^{3m}L\psi_m(x,y,z)$$ for
all $x,y,z \in A$. Moreover, if $f : A\to X$ is a mapping such
that
$$ \|\Delta_m{f(x,y)}\|\leq \theta\|y\|^r,$$ \quad
$$\|f([x,y,z])-[f(x),y^m,z^m]-[x^m,f(y),z^m]-[x^m,y^m,f(z)]\|\leq\psi_m(x,y,z)$$ for all $x,y,z
\in A,$ then $f$ is a ternary $m$-derivation.}
\end{rk}\vskip 5mm
%%%%%%%%%%%%%%%%%%%%%%%%%%%%%%%%%%%%%%%%%%%%%%%%%%%%%%%%%%%%%%%%%%%%%%%%%%%%%%%%%%%%%%%%%%%%%%%%%%%%%%%%%%%%%%%%%%%%%%%
\section {\bf{ Approximation of ternary $m$-$\sigma$-Homomorphism between ternary algebras
}}\vskip 2mm
In this section, we investigate the generalized  stability of
ternary $m$-$\sigma$-Homomorphism between ternary Banach algebras
for the functional equation $(\ref{e12})$.

Throughout this section, we suppose that $A,B$ are two
ternary Banach algebra.
From now on, let $m$ be a positive integer less than $5$ and $\sigma$ a permutation of $\{1,2,3\}.$ \vskip
2mm
\begin{thm}\label{t3.10} Let $f:A\rightarrow B$ be a mapping for which there exist functions $\varphi_m:A\times A \rightarrow  [0,\infty)$
and $\psi_m:A\times A\times A\rightarrow  [0,\infty)$ such that
\begin{equation}\label{e310}
\|\Delta_m{f(x_1,x_2)}\|\leq \varphi_m(x_1,x_2),
\end{equation}
\begin{equation}\label{e320}
\|f([x_1,x_2,x_3])-[f(x_{\sigma(1)}),f(x_{\sigma(2)}),f(x_{\sigma(3)})]\|\leq\psi_m(x_1,x_2,x_3)
\end{equation}
for all $x_1,x_2,x_3 \in A.$ If there exists a constant $0<L<1$ such that
\begin{equation}\label{e330}
\varphi_m\Big(\frac{x_1}{a},\frac{x_2}{a}\Big)\leq
\frac{L}{|a|^m}\varphi_m(x_1,x_2),\end{equation}
\begin{equation}\label{e340}
\psi_m\Big(\frac{x_1}{a},\frac{x_2}{a},\frac{x_3}{a}\Big)\leq
\frac{L}{|a|
^{3m}}\psi_m(x_1,x_2,x_3)
\end{equation}

for all $x_1,x_2,x_3 \in A$ then there exists a unique ternary
$m$-$\sigma$-hommomorphism $H:A\rightarrow B$ such that
\begin{equation}\label{e350}
\|f(x_1)-H(x_1)\| \leq \frac{L}{2|a|^m(1-L)}\varphi_m(x_1,0)
\end{equation}
for all $x_1 \in A.$
\end{thm}\vskip 2mm

\begin{proof} Let us define $\Omega,d$ and $T:\Omega \rightarrow \Omega$ by the same definitions as in the proof of Theorem 3.1, one can show that $T$ has a unique fixed point $H$ in $\Omega$ such that $H(x_1)=\lim_{n\rightarrow\infty}a^{mn}f\Big(\frac{x_1}{a^{n}}\Big)$ that

\begin{equation*}
\|\Delta_m H(x_1,x_2)\|=
\lim_{n\rightarrow\infty}|a|^{mn}\Big\|\Delta_m
f\Big(\frac{x_1}{a^{n}},\frac{x_2}{a^{n}}\Big)\Big\|\leq
\lim_{n\rightarrow\infty}|a|^{mn}\varphi_m\Big(\frac{x_1}{a^{n}},\frac{x_2}{a^{n}}\Big)=0
\end{equation*}
for all $x_1,x_2 \in A$ and so $\Delta_m H(x_1,x_2)=0$.
 By the result in \cite{Es-A-Kh}, $H$ is $m-$mapping.
  On the other hand, it follows from the
definition of $H$ that
\begin{align*}
&\quad\|H([x_1,x_2,x_3])-[H(x_{\sigma(1)}),H(x_{\sigma(2)}),H(x_{\sigma(3)})]\|\\
&=\lim_{n\rightarrow\infty}|a|^{3mn}\Big\|f\Big(\frac{[x_1,x_2,x_3]}{a^{3n}}\Big)-\Big[f\Big(\frac{x_{\sigma(1)}}{a^n}\Big),f\Big(\frac{x_{\sigma(2)}}{a^n}\Big),f\Big(\frac{x_{\sigma(3)}}{a^n}\Big)\Big]\Big\|\\
&
\leq\lim_{n\rightarrow\infty}|a|^{3mn}\psi_m\Big(\frac{x_1}{a^n},\frac{x_2}{a^n},\frac{x_3}{a^n}\Big)=0
\end{align*}
 for all $x_1,x_2,x_3\in A$ and so

  $H([x_1,x_2,x_3])=[H(x_{\sigma(1)}),H(x_{\sigma(2)}),H(x_{\sigma(3)})]$.

  According to Theorem 2.5, since $H$ is the unique fixed
point of $T$ in the set $\Lambda=\{g\in\Omega:d(f,g)<\infty\},$
$H$ is the unique mapping such that
\begin{equation*}
 \|f(x_1)-H(x_1)\|\leq K~\varphi_m(x_1,0)
\end{equation*}
for all $x_1\in A$ and $K>0$. Again, using Theorem 2.5, we have
\begin{equation*}
d(f,H) \leq \frac{1}{1-L}d(f,Tf)\leq
\frac{L}{2|a|^m(1-L)}
\end{equation*}
and so
\begin{equation*}
\|f(x_1)-H(x_1)\| \leq\frac{L}{2|a|^m(1-L)}~\varphi_m(x_1,0)
\end{equation*}
for all $x_1 \in A.$  This completes the proof.
\end{proof}
\vskip 2mm
%%%%%%%%%%%%%%%%%%%%%%%%%%%%%%%%%%%%%%%%%%%%%%%%%%%%%%%%%%%%%%%%%%%%%%%%%%%%%%%%%%%%%%%%%%%%%%%%%%%%%%%%%%%%%%%%%%%%%%%%%%
\begin{cor}\label{c3.20} Let $\theta,r,p$ be non--negative
real numbers with $r,p>m$ and $\frac{3p-r}{2}\geq m$. Suppose that
$f:A\rightarrow B$ is a mapping such that
\begin{equation}\label{b39}
\|\Delta_m f(x_1,x_2)\|\leq \theta(\|x_1\|^r+\|x_2\|^r),
\end{equation}
\begin{equation}\label{b39+}
\|f([x_1,x_2,x_3])-[f(x_{\sigma(1)}),f(x_{\sigma(2)}),f(x_{\sigma(3)})]\|\leq
\theta(\|x_1\|^p.\|x_2\|^p.\|x_3\|^p)
\end{equation}
for all $x_1,x_2,x_3\in A.$ Then there exists a unique ternary
$m$-$\sigma$-homomorphism $H:A\rightarrow B$ satisfying
\begin{equation*} \|f(x_1)-H(x_1)\|
\leq\frac{\theta}{2(|a|^r-|a|^m)}\|x_1\|^r
\end{equation*}
for all $x_1 \in A.$
\end{cor}\vskip 2mm

\begin{proof}
The proof follows from Theorem \ref{t3.10} by taking
\begin{equation*}
\varphi_m(x_1,x_2):=\theta(\|x_1\|^r+\|x_2\|^r),\quad
\psi_m(x_1,x_2,x_3):=\theta(\|x_1\|^p\cdot\|x_2\|^p\cdot\|x_3\|^p)
\end{equation*}
for all $x_1,x_2,x_3\in A$. Then we can choose $L=|a|^{m-r}$ and so the
desired conclusion follows.
\end{proof}\vskip 2mm
%%%%%%%%%%%%%%%%%%%%%%%%%%%%%%%%%%%%%%%%%%%%%%%%%%%%%%%%%%%%%%%%%%%%%%%%%%%%%%%%%%%%%%%%%%%%%%%%%%%%%%%%%%%%%%%%%%%%%%%
\begin{rk}\label{rk}
{\rm Let  $f : A \to B $ be a mapping with $f(0)=0$ such that there
exist functions $\varphi_m:A\times A\rightarrow  [0,\infty)$ and
$\psi_m:A\times A\times A\rightarrow  [0,\infty)$ satisfying
(\ref{e310}) and $(\ref{e320})$. Let $0<L<1$ be a constant such that
$$\varphi_m(ax_1,ax_2)\leq |a|^mL\varphi_m(x_1,x_2),\quad
\psi_m(ax_1,ax_2,ax_3)\leq |a|^{3m}L\psi_m(x_1,x_2,x_3)
$$
for all $x_1,x_2,x_3 \in A.$ By the similar method as in the proof of
Theorem \ref{t3.10}, one can show that there exists a unique
ternary $m$-$\sigma$-homomorphism $H:A\rightarrow B$
satisfying
\[\|f(x_1)-H(x_1)\| \leq
\frac{1}{2|a|^m(1-L)}\varphi_m(x_1,0)\] for all $x_1\in A.$
 For the case
 $$\varphi_m(x_1,x_2):=\theta(\|x_1\|^r+\|x_2\|^r),\quad
 \psi_m(x_1,x_2,x_3):=\theta(\|x_1\|^p+\|x_2\|^p+\|x_3\|^p),
 $$
  where $\theta$ is non--negative real
 numbers
 and $0<r<m$ , $0<p<3m$ and $\frac{p-r}{2}\leq m,$
there exists a unique ternary $m$-$\sigma$-homomorphism
$H:A\rightarrow B$ satisfying
\begin{equation*}
\|f(x_1)-H(x_1)\|
\leq\frac{\theta}{2(|a|^m-|a|^r)}\|x_1\|^r
\end{equation*}
 for all $x_1 \in A.$}
\end{rk}\vskip 2mm
%%%%%%%%%%%%%%%%%%%%%%%%%%%%%%%%%%%%%%%%%%%%%%%%%%%%%%%%%%%%%%%%%%%%%%%%%%%%%%%%%%%%%%%%%%%%%%%%%%%%%%%
 {\small
%----------------------------------------------------------------------%

}
\end{document}